\documentclass[11pt]{amsart}
\usepackage[T1]{fontenc}
\usepackage{amsmath,amssymb,mathtools}
\usepackage[margin=1.22in]{geometry}
\usepackage{hyperref}
\hypersetup{hidelinks}
\newtheorem{theorem}{Theorem}
\newtheorem{lemma}[theorem]{Lemma}
\newcommand{\K}{\mathbb K}
\newcommand{\R}{\mathbb R}
\newcommand{\N}{\mathbb N}
\newcommand{\id}{\operatorname{id}}
\newcommand{\tr}{\operatorname{tr}}
\newcommand{\norm}[1]{\left\|#1\right\|}
\title{Asymptotic projection constants of Lorentz sequence spaces}
\author[D. N\'u\~nez-Alarc\'on]{Daniel N\'u\~nez-Alarc\'on}
\address{Departamento de Matem\'aticas, Facultad de Ciencias, Universidad Nacional de Colombia, Sede Bogot\'a, Bogot\'a, Colombia}
\email{dnuneza@unal.edu.co}
\author[D. M. Pellegrino]{Daniel Marinho Pellegrino}
\address{Departamento de Matem\'atica, Centro de Ci\^encias Exatas e da Natureza, Universidade Federal da Para\'iba, Jo\~ao Pessoa, PB, Brazil}
\email{daniel.pellegrino@academico.ufpb.br}
\author[D. M. Serrano-Rodr\'iguez]{Diana Marcela Serrano-Rodr\'iguez}
\address{Departamento de Matem\'aticas, Facultad de Ciencias, Universidad Nacional de Colombia, Sede Bogot\'a, Bogot\'a, Colombia}
\email{diserranor@unal.edu.co}
\date{}
\begin{document}
\begin{abstract}
We settle a question raised by Defant, Galicer, Mansilla, Masty\l o and Muro (arXiv:2509.26326v2) on the order of growth of the absolute projection constants of finite-dimensional Lorentz sequence spaces $\ell_{2,s}^n$, $1<s<2$. We determine the exact leading coefficients over both scalar fields, including the endpoint $s=1$. Writing $\lambda(X)$ for the absolute projection constant of $X$, we prove that, for each fixed $1\le s<2$,
\[
 \lambda(\ell_{2,s}^n(\mathbb R))\sim\sqrt{\frac{2sn}{\pi(2-s)\log\log n}},\qquad
 \lambda(\ell_{2,s}^n(\mathbb C))\sim\sqrt{\frac{\pi sn}{4(2-s)\log\log n}},
\]
as $n\to\infty$. We also establish the uniform estimate
\[
 \lambda(\ell_{2,s}^n(\mathbb K))\asymp
 \sqrt{\frac{n}{1+(2/s-1)\log(e+\log n)}},
\]
for all $n\ge1$, $1\le s\le2$ and $\mathbb K\in\{\mathbb R,\mathbb C\}$, with comparison constants independent of $n$, $s$ and $\mathbb K$. In particular, the transition to the Euclidean order $\sqrt n$ occurs on the scale $2-s\asymp1/\log\log n$.
\end{abstract}
\maketitle
\section{Introduction}
Let $X$ be a finite-dimensional Banach space over $\K\in\{\R,\mathbb C\}$. The absolute projection constant of $X$ is defined by
\begin{equation}\label{eq:lambda}
 \lambda(X)=\sup_{Y,J}\inf\big\{\norm{P}:P\in\mathcal L(Y,X),\ PJ=\id_X\big\},
\end{equation}
where $J$ runs through the linear isometries of $X$ into Banach spaces $Y$ over $\K$, $\mathcal L(Y,X)$ is the space of bounded linear operators from $Y$ to $X$, and $\id_X$ is the identity on $X$. By the Kadets--Snobar theorem \cite{KSN}, $1\le\lambda(X)\le\sqrt{\dim X}$ for $X\ne\{0\}$. The problem of determining the order of $\lambda(X)$ for particular classes of finite-dimensional spaces is closely related to the theory of operator ideals (see \cite{DJT,TJ}). For spaces with a $1$-symmetric basis, averaging over permutations and Khintchine-type inequalities play an important role; see \cite{KST,KS,KS2}.

The connection between local Banach space theory and the Bohr phenomenon goes back to Defant, Garc\'ia and Maestre \cite{DGM}; see also \cite{BPS} for the asymptotic value of the Bohr radius of the polydisc. Projection constants have also been studied for spaces of Dirichlet polynomials and for the trace class \cite{DGMMMDir,DGMMMTrace}. The question considered here arose in the study of projection, Gordon--Lewis and unconditional basis constants of spaces of analytic polynomials on Lorentz sequence spaces by Defant, Galicer, Mansilla, Masty\l o and Muro \cite{DGMMM}.

For $1\le s\le2$ and $n\in\N$, let $\ell_{2,s}^n$ denote $\K^n$ endowed with the norm
\begin{equation}\label{eq:lorentz}
 \norm{x}_{2,s}=\left(\sum_{i=1}^n i^{s/2-1}(x_i^*)^s\right)^{1/s},
\end{equation}
for $x=(x_i)_{i=1}^n\in\K^n$, where $x^*=(x_i^*)_{i=1}^n$ is the non-increasing rearrangement of $(|x_i|)_{i=1}^n$. Since the weights $i^{s/2-1}$ are non-increasing, the rearrangement inequality gives
\begin{equation}\label{eq:maxnorm}
 \norm{x}_{2,s}=\max_{\sigma\in S_n}\left(\sum_{i=1}^n i^{s/2-1}|x_{\sigma(i)}|^s\right)^{1/s},
\end{equation}
where $S_n$ denotes the symmetric group on $\{1,\ldots,n\}$. In particular, $\norm{\cdot}_{2,s}$ is a norm and the canonical basis is $1$-symmetric.

We write $A\lesssim B$ (or $B\gtrsim A$) if $A\le CB$ for a numerical constant $C$, and $A\asymp B$ if $A\lesssim B$ and $B\lesssim A$. Numerical constants are independent of $n$, $s$ and the scalar field. A subscript $s$ permits dependence on $s$. All logarithms are natural.

For real scalars, Kwapie\'n and Sch\"utt \cite[Theorem~3.2]{KS} proved in 1985 that
\begin{equation}\label{eq:KS}
 \lambda(\ell_{2,1}^n)\asymp \sqrt{\frac{n}{\log(e+\log n)}}.
\end{equation}
Proposition~3.1 of arXiv:2509.26326v2 \cite{DGMMM} gives the order of $\lambda(\ell_{r,s}^n)$ over the complex field for $1<r<\infty$ and $1\le s\le\infty$, except when $r=2$ and $1<s<2$. In particular, $\lambda(\ell_{2,s}^n)\asymp_s\sqrt n$ for $2\le s\le\infty$. In the paragraph following that proposition the authors ask for the order of $\lambda(\ell_{2,s}^n)$ in the remaining range $1<s<2$. Their Lemma~3.2 gives, by comparison with the Euclidean space and with $\ell_{2,1}^n$, the lower bounds
\[
 \lambda(\ell_{2,s}^n)\gtrsim_s
 \begin{cases}
 \displaystyle\frac{1}{(1+\log n)^{1-1/s}}\sqrt{\frac{n}{\log(e+\log n)}},&1<s<4/3,\\[2mm]
 \displaystyle\frac{\sqrt n}{(1+\log n)^{1/s-1/2}},&4/3\le s<2.
 \end{cases}
\]

For each fixed $1<s<2$, the order is the same as in \eqref{eq:KS}. More precisely, the coefficient of $\sqrt{n/\log\log n}$ can be determined exactly for every $1\le s<2$, over both scalar fields.

We write $\ell_{2,s}^n(\K)$ when the scalar field in \eqref{eq:lorentz} is to be specified. For $1\le s<2$ put
\begin{equation}\label{eq:leadingconstant}
 K_s^{\R}=\sqrt{\frac{2s}{\pi(2-s)}},\qquad
 K_s^{\mathbb C}=\sqrt{\frac{\pi s}{4(2-s)}}.
\end{equation}

\begin{theorem}\label{thm:exact}
For every $s\in[1,2)$ and each $\K\in\{\R,\mathbb C\}$, with $s$ fixed as $n\to\infty$,
\begin{equation}\label{eq:normalizedlimit}
 \lim_{n\to\infty}\frac{\lambda(\ell_{2,s}^n(\K))}{K_s^{\K}\sqrt{n/\log\log n}}=1.
\end{equation}
\end{theorem}

At $s=1$, the constants are $\sqrt{2/\pi}$ and $\sqrt\pi/2$ in the real and complex cases, respectively. These are the same constants as in the asymptotic formula for the classical spaces $\ell_p^n$, $1\le p\le2$, obtained by K\"onig, Sch\"utt and Tomczak-Jaegermann \cite[Theorem~1]{KST}: the limit of $\lambda(\ell_p^n)/\sqrt n$ is $\sqrt{2/\pi}$ over $\R$ and $\sqrt\pi/2$ over $\mathbb C$. For each fixed $1\le s<2$,
\[
 \lim_{n\to\infty}\frac{\lambda(\ell_{2,s}^n(\mathbb C))}{\lambda(\ell_{2,s}^n(\R))}=\frac{\pi}{2\sqrt2}.
\]

The constants in \eqref{eq:leadingconstant} tend to infinity as $s\uparrow2$, whereas $\ell_{2,2}^n$ is Euclidean. To compare the projection constants also when $s$ varies with $n$, put
\begin{equation}\label{eq:uniformscale}
 h(n)=\log(e+\log n),\qquad
 G(n,s)=\sqrt{\frac{n}{1+(2/s-1)h(n)}}
 \quad(n\in\N,\ 1\le s\le2).
\end{equation}
Here $h(n)\ge1$ for all $n$, and $h(n)\sim\log\log n$ as $n\to\infty$.

\begin{theorem}\label{thm:main}
There exist constants $0<c\le C<\infty$ such that, for every $n\in\N$, every $s\in[1,2]$ and both $\K\in\{\R,\mathbb C\}$,
\begin{equation}\label{eq:uniformbounds}
 c\sqrt{\frac{n}{1+(2/s-1)h(n)}}
 \le\lambda(\ell_{2,s}^n(\K))
 \le C\sqrt{\frac{n}{1+(2/s-1)h(n)}}.
\end{equation}
\end{theorem}

One may take $c=1/40$ and $C=84$ in \eqref{eq:uniformbounds}.

Defant, Galicer, Mansilla, Masty\l o and Muro \cite{DGMMM} also ask for the optimal dependence of $\lambda(\ell_{r,s}^n)$ on the Lorentz parameters. Theorem~\ref{thm:main} determines this dependence, up to numerical constants, in the case $r=2$ and $1\le s\le2$, including when $s$ varies with the dimension.

For any sequence $(s_n)_{n\ge1}$ in $[1,2]$, write $d_n=(2/s_n-1)h(n)$. Dividing \eqref{eq:uniformbounds} by $\sqrt n$ gives
\begin{equation}\label{eq:transitionbounds}
 \frac{c}{\sqrt{1+d_n}}\le
 \frac{\lambda(\ell_{2,s_n}^n(\K))}{\sqrt n}
 \le\frac{C}{\sqrt{1+d_n}}.
\end{equation}
Thus the projection constants are comparable to $\sqrt n$ if and only if $(d_n)$ is bounded, or equivalently, $2-s_n=O(1/h(n))$. The scale of the transition is therefore $2-s\asymp1/\log\log n$. If $d_n\to\infty$, then
\[
 \lambda(\ell_{2,s_n}^n(\K))\asymp\sqrt{\frac{n}{(2/s_n-1)h(n)}}
\]
for all sufficiently large $n$. For example, $s_n=2-1/h(n)$ gives $d_n=1/s_n\le1$, so the order is $\sqrt n$ even though $s_n<2$ for every $n$. 

\section{Preliminaries}
If $X$ is a finite-dimensional Banach space, the absolutely $1$-summing norm of the identity $\id_X\colon X\to X$ is the least constant $K$ such that
\begin{equation}\label{eq:pi}
 \sum_{j=1}^m\norm{x_j}_X\le K\sup_{\varphi\in B_{X^*}}\sum_{j=1}^m|\varphi(x_j)|
\end{equation}
for all $m\in\N$ and $x_1,\dots,x_m\in X$. We denote this constant by $\pi_1(\id_X)$; as usual, $B_{X^*}$ is the closed unit ball of the dual of $X$.

The classical identity relating $\lambda(X)$ to $\pi_1(\id_X)$ for spaces with a symmetric basis is used by Kwapie\'n and Sch\"utt \cite[Lemma~3.1]{KS}. 
\begin{lemma}\label{lem:identity}
Let $X=(\K^n,\norm{\cdot})$, where $n\in\N$, be a Banach space whose norm is invariant under coordinate permutations and sign changes. Then
\begin{equation}\label{eq:identity}
 \lambda(X)\pi_1(\id_X)=n.
\end{equation}
\end{lemma}
\begin{proof}
Put $K=B_{X^*}$ with its norm topology and let $J\colon X\to C(K)$ be the isometry given by $Jx(\varphi)=\varphi(x)$. Here $C(K)$ is the space of continuous $\K$-valued functions on $K$, with the supremum norm. Let $P\colon C(K)\to X$ be bounded and satisfy $PJ=\id_X$. Such a map exists by Hahn--Banach applied to the coordinate functionals on $J(X)$. Given $\delta>0$, choose a finite cover of $K$ by norm balls of radius $\delta$ centered at $\varphi_1,\dots,\varphi_M\in K$, and a continuous nonnegative partition of unity $h_1,\dots,h_M$ subordinate to this cover. The map $T_\delta\colon X\to C(K)$ given by
\[
 T_\delta x=\sum_{j=1}^M\varphi_j(x)h_j
\]
satisfies $\norm{T_\delta-J}\le\delta$, since $h_j(\varphi)>0$ only if $\norm{\varphi-\varphi_j}<\delta$ and $\sum_j h_j(\varphi)=1$. Hence $\norm{PT_\delta-\id_X}\le\delta\norm{P}$. Write $v_j=P h_j$. For $\psi\in B_{X^*}$, choose scalars $\theta_j$ of modulus one so that $\theta_j\psi(v_j)=|\psi(v_j)|$. Since $\sum_j h_j=1$,
\begin{equation*}
 \sum_{j=1}^M|\psi(v_j)|
 =\left|\psi\left(P\sum_{j=1}^M\theta_j h_j\right)\right|
 \le\norm{P}.
\end{equation*}
By \eqref{eq:pi} and the trace formula for rank-one operators,
\[
 |\tr(PT_\delta)|\le\sum_{j=1}^M|\varphi_j(v_j)|
 \le\sum_{j=1}^M\norm{v_j}\le\pi_1(\id_X)\norm{P}.
\]
Since the trace is continuous on $\mathcal L(X,X)$, letting $\delta\to0$ gives $n\le\pi_1(\id_X)\norm{P}$. Taking the infimum over $P$ and using \eqref{eq:lambda} gives $n\le\pi_1(\id_X)\lambda(X)$.

For the reverse inequality let $G$ be the finite group of signed coordinate permutations, acting isometrically on $X$. Fix a finite family $x_1,\dots,x_m$ with $S=\sum_j\norm{x_j}>0$ and choose $\varphi_j\in B_{X^*}$ with $\varphi_j(x_j)=\norm{x_j}$. For any isometric embedding $J\colon X\to Y$ and every $g\in G$, extend the functional $\varphi_j\circ g^{-1}\circ J^{-1}$ on $J(X)$ to a functional $f_{j,g}\in B_{Y^*}$ by Hahn--Banach. Define $Q\colon Y\to X$ by
\[
 Qy=\frac{n}{S|G|}\sum_{g\in G}\sum_{j=1}^m f_{j,g}(y)gx_j.
\]
If $T\colon X\to X$ is given by $Tx=\sum_{j=1}^m\varphi_j(x)x_j$, then
\[
 QJ=\frac{n}{S|G|}\sum_{g\in G}gTg^{-1}.
\]
The average commutes with $G$. Commutation with the coordinate sign changes makes its matrix diagonal, and commutation with the permutations makes its diagonal entries equal. Since $\tr(T)=S$, we have $\tr(QJ)=n$.
Hence $QJ=\id_X$. Set $D=\sup_{\psi\in B_{X^*}}\sum_j|\psi(x_j)|$. For $y\in B_Y$ and $\psi\in B_{X^*}$ we have
\[
 |\psi(Qy)|\le\frac{n}{S|G|}\sum_{g\in G}\sum_{j=1}^m|\psi(gx_j)|\le\frac{nD}{S},
\]
since $\psi\circ g\in B_{X^*}$. Thus $\lambda(X)\le nD/S$. Taking families with $S/D$ arbitrarily close to $\pi_1(\id_X)$ proves \eqref{eq:identity}.
\end{proof}

We identify the dual of $\ell_{2,s}^n$ with $\K^n$ by the bilinear pairing $\langle x,y\rangle=\sum_{i=1}^n x_i y_i$. For $y\in\K^n$ put
\[
 \norm{y}_*=\sup_{\norm{x}_{2,s}\le1}|\langle x,y\rangle|.
\]

For $y\in B_{X^*}$, with $X=\ell_{2,s}^n$ and $1\le s\le2$, testing against the vector of the first $j$ ones gives
\begin{equation}\label{eq:ycoord}
 j y_j^*\le\left(\sum_{h=1}^j h^{s/2-1}\right)^{1/s}\le (2/s)^{1/s}\sqrt j.
\end{equation}
The last inequality follows by integration of $r^{s/2-1}$ over $(0,j]$.

\section{Averages over permutations}
For the symmetric group $S_n$, let $\mu_n$ be the probability measure assigning mass $1/n!$ to each permutation. For $x,y\in\K^n$ put
\begin{equation}\label{eq:F}
 F_n(x,y)=\int_{S_n}\left(\sum_{i=1}^n|x_i y_{\sigma(i)}|^2\right)^{1/2}\,d\mu_n(\sigma).
\end{equation}
The map $F_n\colon\K^n\times\K^n\to[0,\infty)$ is invariant under separate permutations of the absolute values of the coordinates of $x$ and $y$. Estimates for these averages are central to the work of Kwapie\'n and Sch\"utt \cite[Theorem~1.2]{KS}; see also \cite{KS2}. Let $\nu_n$ be the uniform probability measure on the group $G$ of signed coordinate permutations, consisting of the maps $g\colon\K^n\to\K^n$ of the form $(gx)_{\sigma(i)}=\varepsilon_i x_i$ for $\sigma\in S_n$ and $\varepsilon\in\{-1,1\}^n$. We have
\[
 \int_G|\langle gx,y\rangle|\,d\nu_n(g)
 =\frac{1}{2^n n!}\sum_{\sigma\in S_n}\sum_{\varepsilon\in\{-1,1\}^n}
 \left|\sum_{i=1}^n\varepsilon_i x_i y_{\sigma(i)}\right|.
\]
By the Khintchine inequality,
\begin{equation}\label{eq:Khinchine}
 c F_n(x,y)\le\int_G|\langle gx,y\rangle|\,d\nu_n(g)\le F_n(x,y),
\end{equation}
where $c=1/\sqrt3$ works over both scalar fields. Indeed, let $S=\sum_{i=1}^n\xi_i a_i$, where the $\xi_i$ are independent uniform signs, or, over $\mathbb C$, independent uniform variables on the unit circle. Then $\mathbb E|S|^2=\sum_i|a_i|^2$ and
\[
 \mathbb E|S|^4\le\sum_i|a_i|^4+6\sum_{i<j}|a_i|^2|a_j|^2\le3(\mathbb E|S|^2)^2.
\]
In the fourth moment, terms with an index occurring only once vanish. By H\"older,
\[
 \mathbb E|S|^2\le(\mathbb E|S|)^{2/3}(\mathbb E|S|^4)^{1/3}.
\]
Together with Cauchy--Schwarz this gives, for either choice of the $\xi_i$,
\begin{equation}\label{eq:phasekhintchine}
 \frac1{\sqrt3}\left(\sum_i|a_i|^2\right)^{1/2}\le\mathbb E|S|\le\left(\sum_i|a_i|^2\right)^{1/2}.
\end{equation}

Let $X=\ell_{2,s}^n$. If $a>0$, $y\in B_{X^*}$ and $F_n(x,y)\ge a\norm{x}_X$ for all $x\in X$, then for any finite family $x_1,\dots,x_m$ in $X$,
\begin{align*}
 \sup_{z\in B_{X^*}}\sum_{j=1}^m|\langle x_j,z\rangle|
 &\ge\int_G\sum_{j=1}^m|\langle gx_j,y\rangle|\,d\nu_n(g)\\
 &\ge ca\sum_{j=1}^m\norm{x_j}_X.
\end{align*}
The first inequality uses invariance of the dual unit ball under the action induced by $G$. Consequently
\begin{equation}\label{eq:upperpi}
 \pi_1(\id_X)\le (ca)^{-1}.
\end{equation}
Conversely, using the finite family $(gx)_{g\in G}$ for a fixed $x\ne0$ in \eqref{eq:pi}, we obtain
\[
 \norm{x}_X\le\pi_1(\id_X)\sup_{y\in B_{X^*}}\int_G|\langle gx,y\rangle|\,d\nu_n(g)
 \le\pi_1(\id_X)\sup_{y\in B_{X^*}}F_n(x,y).
\]
Thus
\begin{equation}\label{eq:lowerpi}
 \pi_1(\id_X)\ge\frac{\norm{x}_X}{\sup_{y\in B_{X^*}}F_n(x,y)}.
\end{equation}

For a nonnegative matrix $W=(w_{ij})_{i,j=1}^n$, put $Z\colon S_n\to[0,\infty)$, $Z(\sigma)=\sum_{i=1}^n w_{i,\sigma(i)}$, and $M=\int_{S_n}Z\,d\mu_n$.
\begin{lemma}\label{lem:moment}
If $n\ge2$ and $\max_{i,j}w_{ij}\le2M$, then
\[
 \int_{S_n}\sqrt{Z}\,d\mu_n\ge\frac12\sqrt{M}.
\]
\end{lemma}
\begin{proof}
For $1\le i,j\le n$, there are $(n-1)!$ permutations satisfying $\sigma(i)=j$. For $i\ne i'$ and $j\ne j'$, there are $(n-2)!$ satisfying $\sigma(i)=j$ and $\sigma(i')=j'$. Thus
\begin{equation}\label{eq:momentmean}
 M=\sum_{i,j=1}^n w_{ij}\frac{(n-1)!}{n!}=\frac1n\sum_{i,j=1}^n w_{ij}.
\end{equation}
Moreover,
\begin{align*}
 \int_{S_n}Z^2\,d\mu_n
 &=\frac1n\sum_{i,j}w_{ij}^2+\frac1{n(n-1)}\sum_{i\ne i'}\sum_{j\ne j'}w_{ij}w_{i'j'}\\
 &\le(\max_{i,j}w_{ij})M+\frac n{n-1}M^2.
\end{align*}
In particular, for any nonnegative matrix $W$,
\begin{equation}\label{eq:generalvariance}
 \int_{S_n}(Z-M)^2\,d\mu_n\le (\max_{i,j} w_{ij})M+\frac{M^2}{n-1}.
\end{equation}
Under the hypothesis of the lemma, $\int_{S_n} Z^2\,d\mu_n\le4M^2$. If $M=0$, the assertion holds. For $M>0$, H\"older's inequality with exponents $3/2$ and $3$ implies
\[
 M=\int_{S_n}(\sqrt Z)^{2/3}(Z^2)^{1/3}\,d\mu_n
 \le\left(\int_{S_n}\sqrt Z\,d\mu_n\right)^{2/3}
     \left(\int_{S_n}Z^2\,d\mu_n\right)^{1/3}.
\]
Consequently
\[
 \int_{S_n}\sqrt Z\,d\mu_n\ge\frac{M^{3/2}}{(\int_{S_n}Z^2\,d\mu_n)^{1/2}}\ge\frac12\sqrt M.
\]
\end{proof}

\section{The asymptotic constant}
Fix $1\le s<2$ and let $t$ be its conjugate exponent, so that $t=s/(s-1)$ if $s>1$ and $t=\infty$ if $s=1$. Put $\beta=2/s-1=1-2/t$, with the convention $1/\infty=0$. Throughout this section, $X=\ell_{2,s}^n$ and $n\to\infty$ with $s$ fixed. The implicit constants and the rates of convergence may depend on $s$; the estimates are uniform on the indicated unit balls.

Let $\mathbb T_{\K}=\{z\in\K:|z|=1\}$, and let $\rho$ be its normalized Haar measure. Set $\Omega_n=\mathbb T_{\K}^n$ and $\rho_n=\rho^{\otimes n}$. For $x,y\in\K^n$ define
\begin{equation}\label{eq:phaseaverage}
 A_n(x,y)=\int_{S_n}\int_{\Omega_n}
 \left|\sum_{i=1}^n\omega_i x_i y_{\sigma(i)}\right|\,d\rho_n(\omega)\,d\mu_n(\sigma).
\end{equation}
This average is invariant under rearrangements of the absolute values of either vector. For a Borel probability measure $\eta$ on $B_{X^*}$ and $a>0$, if
\begin{equation}\label{eq:common}
 \int_{B_{X^*}} A_n(x,y)\,d\eta(y)\ge a\norm{x}_{2,s}\qquad(x\in X),
\end{equation}
then $\pi_1(\id_X)^{-1}\ge a$. Indeed, for $x_1,\ldots,x_m\in X$, the supremum in \eqref{eq:pi} is at least the integral of the sum over all $y$, permutations and phases in \eqref{eq:common}, since each corresponding functional belongs to $B_{X^*}$. On the other hand, for each $x\ne0$,
\begin{equation}\label{eq:orbitexact}
 \frac{1}{\pi_1(\id_X)}\le \frac{1}{\norm{x}_{2,s}}\sup_{y\in B_{X^*}}A_n(x,y).
\end{equation}
In the real case, apply \eqref{eq:pi} to the finite orbit of $x$ under signed permutations, counted with multiplicity. In the complex case, do the same with permutations and the $m$th roots of unity in each coordinate. The integrals over these roots converge to the integral in \eqref{eq:phaseaverage} uniformly in $y\in B_{X^*}$: replacing each phase by a nearest root changes the integrand by at most $(2\pi/m)\sum_i|x_i y_{\sigma(i)}|\le (2\pi/m)\norm{x}_{2,s}$. Letting $m\to\infty$ proves \eqref{eq:orbitexact}.

Let $g$ be a real standard Gaussian if $\K=\R$, and $g=(g_1+\mathrm i g_2)/\sqrt2$ for independent real standard Gaussians $g_1,g_2$ if $\K=\mathbb C$. In both cases $\int |g|^2=1$, where integration is with respect to the law of $g$. Denote the mean of $|g|$ by $\gamma_{\K}$. The corresponding integrals are
\begin{equation}\label{eq:gammas}
 \gamma_{\R}=\sqrt{\frac2\pi}\int_0^\infty r e^{-r^2/2}\,dr=\sqrt{\frac2\pi},\qquad
 \gamma_{\mathbb C}=\int_0^\infty 2r^2 e^{-r^2}\,dr=\frac{\sqrt\pi}{2}.
\end{equation}

\begin{lemma}\label{lem:gaussian}
For each $\varepsilon>0$ there is a constant $C_\varepsilon$ such that, for every $n$ and $a_1,\ldots,a_n\in\K$,
\begin{equation}\label{eq:gaussian}
 \left|\int_{\Omega_n}\left|\sum_{i=1}^n\omega_i a_i\right|d\rho_n(\omega)
 -\gamma_{\K}\left(\sum_{i=1}^n|a_i|^2\right)^{1/2}\right|
 \le\varepsilon\left(\sum_{i=1}^n|a_i|^2\right)^{1/2}+C_\varepsilon\max_i|a_i|.
\end{equation}
\end{lemma}
\begin{proof}
Theorem~2 of K\"onig, Sch\"utt and Tomczak-Jaegermann \cite{KST}, for Rademacher variables over $\R$ and Steinhaus variables over $\mathbb C$, gives the stronger estimate
\[
 \left|\int_{\Omega_n}\left|\sum_{i=1}^n\omega_i a_i\right|d\rho_n(\omega)
 -\gamma_{\K}\left(\sum_{i=1}^n|a_i|^2\right)^{1/2}\right|
 \le (1-\gamma_{\K})\max_i|a_i|.
\]
In particular, \eqref{eq:gaussian} holds with $C_\varepsilon=1-\gamma_{\K}$ for each $\varepsilon>0$.
\end{proof}

For all sufficiently large $n$, set
\begin{equation}\label{eq:scales}
 H=\log\log n,\qquad R=10\log H,\qquad Q=R^2,\qquad u_j=\log(n/j)\quad(1\le j\le n).
\end{equation}
As $H\to\infty$, we have $R/Q=1/(10\log H)\to0$, $\log Q=2\log(10\log H)=o(H)$ and $Q^2=10^4(\log H)^4=o(e^H)=o(\log n)$. Define
\begin{equation}\label{eq:weights}
 w_j=(Q+u_j)^{-\beta},\qquad D=\sum_{j=1}^n\frac{1}{j(Q+u_j)}.
\end{equation}
The function $r\mapsto 1/[r(Q+\log(n/r))]$ is decreasing for $Q>1$. Integral comparison and the substitution $u=\log(n/r)$ give
\begin{equation}\label{eq:D}
 D=\int_0^{\log n}\frac{du}{Q+u}+O(Q^{-1})
 =\log\frac{Q+\log n}{Q}+O(Q^{-1})=(1+o(1))H.
\end{equation}

\begin{lemma}\label{lem:weighteddual}
For $Q>2$,
\begin{equation}\label{eq:weighteddual}
 D^\beta\le\sup_{y\in B_{X^*}}\sum_{j=1}^n w_j(y_j^*)^2
 \le(1+4\beta/Q)^2D^\beta.
\end{equation}
The lower bound is realized by the non-increasing vector $y\in B_{X^*}$ given by
\begin{equation}\label{eq:q}
 q_j=\frac{1}{\sqrt j(Q+u_j)^{1/t}},\qquad y=D^{-1/t}q\qquad(1\le j\le n).
\end{equation}
\end{lemma}
\begin{proof}
The function $r\mapsto r^{-1/2}(Q+\log(n/r))^{-1/t}$ is decreasing on $[1,n]$: its logarithmic derivative is $r^{-1}(-1/2+1/[t(Q+\log(n/r))])<0$. For $1<s<2$, rearrangement and H\"older's inequality give
\[
 |\langle x,q\rangle|\le\sum_j x_j^*q_j
 \le\norm{x}_{2,s}\left(\sum_j j^{t(1/s-1/2)}q_j^t\right)^{1/t}
 =D^{1/t}\norm{x}_{2,s}.
\]
For $s=1$, $q_j=j^{-1/2}$ and the same estimate follows directly from the definition of the norm. Thus $y\in B_{X^*}$. In fact, $\norm{y}_*=1$: the decreasing vector with coordinates $j^{-1/2}(Q+u_j)^{-1/s}$ has norm $D^{1/s}$ and pairs with $y$ to give $D^{1/s}$. Since $\beta+2/t=1$,
\[
 \sum_{j=1}^n w_j y_j^2=D^{-2/t}\sum_{j=1}^n\frac{1}{j(Q+u_j)}=D^\beta.
\]

For the upper bound, let $y\in B_{X^*}$ be nonnegative and decreasing. Put $v_j=w_j y_j$, $z_j=\max_{j\le k\le n}v_k$ and $B=\sum_j w_j y_j^2$. For $h\ge2$, monotonicity of $y$ gives
\[
 v_h-v_{h-1}\le (w_h-w_{h-1})y_h
 \le \frac{2\beta}{Qh}v_h.
\]
Indeed, if $w(r)=(Q+\log(n/r))^{-\beta}$, then $w'(r)=\beta w(r)/[r(Q+\log(n/r))]$. For $h-1\le r\le h$, one has $w(r)\le w_h$, $r\ge h/2$ and $Q+\log(n/r)\ge Q$. Integrating $w'$ from $h-1$ to $h$ gives $w_h-w_{h-1}\le 2\beta w_h/(Qh)$. It follows by summation that
\begin{equation}\label{eq:majorant}
 z_j\le v_j+\frac{2\beta}{Q}\sum_{h>j}\frac{v_h}{h}.
\end{equation}
For a nonnegative vector $v$ denote by $\norm{v}_{s,0}=(\sum_{j=1}^n j^{s/2-1}v_j^s)^{1/s}$ the weighted norm without rearrangement. If $b_h=\sqrt h v_h$, then
\[
 \sqrt j\sum_{h>j}\frac{v_h}{h}=\sum_{h>j}\frac1h\sqrt{j/h}\,b_h.
\]
For each $j$, $\sum_{h>j}\sqrt{j/h}/h\le2$, and for each $h$, $\sum_{j<h}\sqrt{j/h}/j\le2$, by comparison with the integrals of $r^{-3/2}$ and $r^{-1/2}$, respectively. H\"older's inequality gives
\begin{align*}
 \sum_j\frac1j\left(\sum_{h>j}\frac{\sqrt{j/h}}h b_h\right)^s
 &\le 2^{s-1}\sum_j\frac1j\sum_{h>j}\frac{\sqrt{j/h}}h b_h^s\\
 &\le 2^s\sum_h\frac{b_h^s}{h}.
\end{align*}
Since $z$ is non-increasing, \eqref{eq:majorant} and Minkowski's inequality imply
\[
 \norm{z}_{2,s}=\norm{z}_{s,0}\le(1+4\beta/Q)\norm{v}_{s,0}.
\]
Also $B=\sum_j y_jv_j\le\sum_j y_jz_j\le\norm{z}_{2,s}$. Applying H\"older's inequality with exponents $2/s$ and $2/(2-s)$, and using $\beta s/(2-s)=1$, we obtain
\begin{align*}
 \norm{v}_{s,0}^s
 &=\sum_j(w_j y_j^2)^{s/2}j^{s/2-1}w_j^{s/2}\\
 &\le B^{s/2}\left(\sum_j\frac{w_j^{s/(2-s)}}j\right)^{1-s/2}
 =B^{s/2}D^{1-s/2}.
\end{align*}
Taking the $s$th root and using $(1-s/2)/s=\beta/2$, we obtain
\[
 B\le\norm{z}_{2,s}\le(1+4\beta/Q)B^{1/2}D^{\beta/2}.
\]
If $B>0$, division by $B^{1/2}$ followed by squaring gives \eqref{eq:weighteddual}; for $B=0$ the estimate is immediate.
\end{proof}

\subsection*{An upper estimate for the permutation averages}
For $s=1$, the vector used by Kwapie\'n and Sch\"utt in the proof of \cite[Theorem~3.2]{KS} has coordinates $1/(\sqrt i\log(2i))$. For $1\le s<2$, consider the vector $p\in X$ given by
\begin{equation}\label{eq:exactp}
 p_i=\frac{1}{\sqrt i\,(Q+\log i)^{1/s}}\qquad(1\le i\le n).
\end{equation}
As $p$ is decreasing, integral comparison gives
\begin{equation}\label{eq:exactpnorm}
 \norm{p}_{2,s}^s=\sum_{i=1}^n\frac{1}{i(Q+\log i)}=\log\frac{Q+\log n}{Q}+O(Q^{-1})=(1+o(1))H.
\end{equation}
For $E\subset\{1,\ldots,n\}^2$ and $x,y\in\K^n$, define
\[
 A_E(x,y)=\int_{S_n}\int_{\Omega_n}\left|\sum_{i:(i,\sigma(i))\in E}\omega_i x_i y_{\sigma(i)}\right|d\rho_n\,d\mu_n.
\]
Put
\begin{align*}
 E_+&=\{(i,j):ij\ge ne^R\},\\
 E_0&=\{(i,j):ne^{-R}\le ij<ne^R\},\\
 E_-&=\{(i,j):ij<ne^{-R}\}.
\end{align*}
By the triangle inequality, $A_n(p,y)\le A_{E_+}(p,y)+A_{E_0}(p,y)+A_{E_-}(p,y)$. Fix $y\in B_{X^*}$, nonnegative and non-increasing. For each $j$, put $a=ne^R/j=e^{u_j+R}$. The sum of $p_i^2$ over $i\ge a$ is bounded by the value of the decreasing function $f(r)=1/[r(Q+\log r)^{2/s}]$ at $a$ plus its integral from $a$ to infinity. Since $2/s=1+\beta$, this gives
\begin{align*}
 \sum_{i:(i,j)\in E_+}p_i^2
 &\le \frac{1}{a(Q+\log a)^{2/s}}+\frac{(Q+\log a)^{-\beta}}{\beta}\\
 &\le \frac{1+\beta/Q}{\beta}(Q+u_j)^{-\beta}.
\end{align*}
In particular, by the permutation counting formula \eqref{eq:momentmean} and Lemma~\ref{lem:weighteddual},
\begin{equation}\label{eq:uppermomentexact}
 \int_{S_n}\sum_{i:(i,\sigma(i))\in E_+}p_i^2 y_{\sigma(i)}^2\,d\mu_n
 \le \frac{1+o(1)}{\beta n}H^\beta.
\end{equation}
For the terms in $E_+$, \eqref{eq:exactp} and \eqref{eq:ycoord} give $|p_i y_j|\le 2(ne^R)^{-1/2}$. Apply Lemma~\ref{lem:gaussian} for each permutation and Jensen's inequality to the square root of the sum in \eqref{eq:uppermomentexact}. For each fixed $\varepsilon>0$,
\begin{equation}\label{eq:upperdense}
 A_{E_+}(p,y)\le (\gamma_{\K}+\varepsilon)\sqrt{\frac{(1+o(1))H^\beta}{\beta n}}+\frac{2C_\varepsilon e^{-R/2}}{\sqrt n}.
\end{equation}

For $E_0$, the possible $i$ for a fixed $j$ are between $a_0=\max(1,e^{u_j-R})$ and $b_0=e^{u_j+R}$. Using the same decreasing function $f$ and the fact that $\log(b_0/a_0)\le2R$, we have
\begin{align*}
 \sum_{i:(i,j)\in E_0}p_i^2
 &\le f(a_0)+\int_{a_0}^{b_0} f(r)\,dr\\
 &\le (1+2R)(Q+\log a_0)^{-2/s}\\
 &\le (1+2R)(1-R/Q)^{-2/s}(Q+u_j)^{-2/s}.
\end{align*}
By Cauchy--Schwarz in the phases and Jensen in the permutations, and by \eqref{eq:weighteddual},
\begin{align}\label{eq:transition}
 A_{E_0}(p,y)^2
 &\le\frac1n\sum_j y_j^2\sum_{i:(i,j)\in E_0}p_i^2\notag\\
 &\le\frac{(1+2R)(1-R/Q)^{-2/s}}{nQ}\sum_j w_j y_j^2\notag\\
 &\le\frac{(1+2R)(1-R/Q)^{-2/s}(1+4\beta/Q)^2D^\beta}{nQ}
 =o(H^\beta/n).
\end{align}

For $E_-$, the triangle inequality gives
\begin{equation}\label{eq:sparsefirst}
 A_{E_-}(p,y)\le\frac1n\sum_{j\le ne^{-R}} y_j\sum_{i\le ne^{-R}/j}p_i.
\end{equation}
For real $m\ge1$ and $Q\ge4/s$, the function $b\colon[1,\infty)\to(0,\infty)$, $b(r)=\sqrt r(Q+\log r)^{-1/s}$, satisfies
\[
 b'(r)=\frac{(Q+\log r)^{-1/s}}{\sqrt r}\left(\frac12-\frac{1}{s(Q+\log r)}\right)\ge\frac{(Q+\log r)^{-1/s}}{4\sqrt r}.
\] The summand in \eqref{eq:exactp} is decreasing, so
\[
 \sum_{i\le m}p_i\le p_1+\int_1^m\frac{dr}{\sqrt r(Q+\log r)^{1/s}}\le 5\sqrt m(Q+\log m)^{-1/s}.
\]
Define $z\in X$ by $z_j=j^{-1/2}(Q+u_j-R)^{-1/s}$ for $j\le ne^{-R}$ and $z_j=0$ otherwise. For $j$ in its support the defining function has logarithmic derivative $(-1/2+1/[s(Q+\log(n/j)-R)])/j<0$. Thus $z$ is non-increasing and, by integral comparison,
\[
 \norm{z}_{2,s}^s=\sum_{j\le ne^{-R}}\frac{1}{j(Q+u_j-R)}
 \le Q^{-1}+\log\frac{Q+\log n-R}{Q}=O(H).
\]
Combining this estimate with \eqref{eq:sparsefirst} and duality gives
\begin{equation}\label{eq:sparse}
 A_{E_-}(p,y)\le\frac{5e^{-R/2}}{\sqrt n}\sum_j y_j z_j\lesssim_s \frac{e^{-R/2}H^{1/s}}{\sqrt n}=o(H^{\beta/2}/\sqrt n).
\end{equation}
The last equality follows from $1/s-\beta/2=1/2$ and $e^{-R/2}=H^{-5}$. By \eqref{eq:exactpnorm}, $\norm{p}_{2,s}=(1+o(1))H^{1/s}$. Since $\beta/2+1/2=1/s$, the estimates \eqref{eq:upperdense}, \eqref{eq:transition} and \eqref{eq:sparse} give, for every fixed $\varepsilon>0$,
\[
 \frac{\sqrt{\beta nH}}{\norm{p}_{2,s}}\sup_{y\in B_{X^*}}A_n(p,y)
 \le\gamma_{\K}+\varepsilon+o(1).
\]
Taking the upper limit and then the infimum over $\varepsilon>0$ gives
\begin{equation}\label{eq:upperorbitlimit}
 \limsup_{n\to\infty}\sqrt{\beta n H}\,\frac{\sup_{y\in B_{X^*}}A_n(p,y)}{\norm{p}_{2,s}}\le\gamma_{\K}.
\end{equation}

\subsection*{A lower estimate for the permutation averages}
Let $y=D^{-1/t}q$ be the functional in Lemma~\ref{lem:weighteddual}, given by \eqref{eq:q}. Put $m=\lceil e^{Q^2}\rceil$, which is less than $n$ for sufficiently large $n$, and $\delta=H^{-1/8}$. Fix a nonnegative non-increasing $x$ with $\norm{x}_{2,s}=1$. Write
\begin{equation}\label{eq:headtail}
 U=\left(\sum_{i<m}i^{s/2-1}x_i^s\right)^{1/s},\qquad
 V=\left(\sum_{i=m}^n i^{s/2-1}x_i^s\right)^{1/s}.
\end{equation}
Then $U^s+V^s=1$.

Suppose first that $U\le\delta$. For each $i\ge m$, put $a=ne^R/i$, so that $e^R\le a\le n$. The function $r\mapsto 1/[r(Q+\log(n/r))^{2/t}]$ is decreasing on $[1,n]$. The sum of its values at the integers $j\ge a$ is at least its integral from $a$ to $n$ minus its value at $a$ (use the integral from $\lceil a\rceil$ to $n$). Substituting $u=\log(n/r)$ gives
\begin{align*}
 \sum_{j:(i,j)\in E_+}q_j^2
 &\ge \int_0^{\log i-R}(Q+u)^{-2/t}\,du- e^{-R}Q^{-2/t}\\
 &=\frac{(Q+\log i-R)^\beta-Q^\beta}{\beta}-e^{-R}Q^{-2/t}\\
 &\ge (1-o(1))\frac{(\log i)^\beta}{\beta}.
\end{align*}
Since $i\ge m$, $\log i\ge Q^2$. The lower bound for the sum, divided by $(\log i)^\beta/\beta$, is at least
\[
 (1-R/Q^2)^\beta-Q^{-\beta}-\beta e^{-R}Q^{-2/t-2\beta},
\]
which tends to $1$ for the fixed $\beta>0$, independently of $i$. Put
\[
 Z(\sigma)=\sum_{i:(i,\sigma(i))\in E_+}x_i^2y_{\sigma(i)}^2,\qquad M=\int_{S_n}Z\,d\mu_n.
\]
By \eqref{eq:momentmean}, \eqref{eq:D} and \eqref{eq:q},
\begin{equation}\label{eq:lowerM}
 M\ge\frac{1-o(1)}{\beta n H^{2/t}}\sum_{i=m}^n x_i^2(\log i)^\beta.
\end{equation}
H\"older's inequality with exponents $2/s$ and $2/(2-s)$ gives
\begin{align*}
 V^s&=\sum_{i=m}^n (x_i^2(\log i)^\beta)^{s/2}i^{s/2-1}(\log i)^{-\beta s/2}\\
 &\le\left(\sum_{i=m}^n x_i^2(\log i)^\beta\right)^{s/2}
 \left(\sum_{i=m}^n\frac{1}{i\log i}\right)^{1-s/2}.
\end{align*}
The last sum is at most $1/(m\log m)+\log\log n-\log\log m\le H+o(1)$. Since $V^s\ge1-\delta^s$, taking the power $2/s$ in the preceding inequality gives
\[
 \sum_{i=m}^n x_i^2(\log i)^\beta
 \ge\frac{(1-\delta^s)^{2/s}}{(H+o(1))^\beta}=(1-o(1))H^{-\beta}.
\]
Substitution in \eqref{eq:lowerM}, with $\beta+2/t=1$, gives
\begin{equation}\label{eq:Mmin}
 M\ge\frac{1-o(1)}{\beta n H}.
\end{equation}
For a nonnegative decreasing $x$, $\sqrt i x_i\le\norm{x}_2$, and
\[
 \norm{x}_{2,s}^s=\sum_{i:x_i>0}x_i^2(\sqrt i x_i)^{s-2}
 \ge\norm{x}_2^{s-2}\sum_i x_i^2=\norm{x}_2^s
\]
when $x\ne0$. Thus $x_i\le1/\sqrt i$. By \eqref{eq:q}, $y_j\le1/\sqrt j$ for large $n$. In particular, $x_i^2y_j^2\le e^{-R}/n$ on $E_+$. Applying \eqref{eq:generalvariance} to $w_{ij}=x_i^2y_j^2\mathbf 1_{E_+}(i,j)$, we obtain
\begin{equation}\label{eq:varbound}
 \int_{S_n}(Z-M)^2\,d\mu_n\le \frac{e^{-R}}n M+\frac{M^2}{n-1}.
\end{equation}
For sufficiently large $n$, \eqref{eq:Mmin} implies $M\ge1/(2\beta nH)$. Hence
\begin{equation}\label{eq:relativevariance}
 \frac{\int_{S_n}(Z-M)^2\,d\mu_n}{M^2}\le 2\beta H e^{-R}+\frac1{n-1}=o(1).
\end{equation}
Using $|\sqrt{Z}-\sqrt M|\le |Z-M|/\sqrt M$ and Cauchy--Schwarz gives
\[
 \left|\int_{S_n}\sqrt Z\,d\mu_n-\sqrt M\right|\le\frac{(\int_{S_n}(Z-M)^2\,d\mu_n)^{1/2}}{\sqrt M}.
\]
By \eqref{eq:relativevariance},
\begin{equation}\label{eq:sqrtconc}
 \int_{S_n}\sqrt Z\,d\mu_n\ge(1-o(1))\sqrt M.
\end{equation}
For each permutation, the integral of the absolute value of the full phase sum is at least that of the sum over the pairs in $E_+$. Indeed, the omitted phases have mean zero, and Jensen's inequality applies conditionally on the retained phases. Apply \eqref{eq:gaussian} to the retained coefficients, whose maximum is at most $e^{-R/2}/\sqrt n$, and use \eqref{eq:sqrtconc}. For any $0<\varepsilon<\gamma_{\K}$,
\begin{align}\label{eq:lowerdense}
 A_n(x,y)&\ge A_{E_+}(x,y)\notag\\
 &\ge (\gamma_{\K}-\varepsilon)(1-o(1))\sqrt M-C_\varepsilon e^{-R/2}/\sqrt n\notag\\
 &\ge\frac{(\gamma_{\K}-\varepsilon)(1-o(1))}{\sqrt{\beta nH}}-\frac{C_\varepsilon e^{-R/2}}{\sqrt n}.
\end{align}

For vectors with $U>\delta$, let $y^{(0)}=n^{-1/2}(1,\ldots,1)$. By Cauchy--Schwarz and $\norm{x}_2\le\norm{x}_{2,s}$, this vector is in $B_{X^*}$. The phase average satisfies
\begin{equation}\label{eq:y0}
 A_n(x,y^{(0)})\ge c_0\frac{\norm{x}_2}{\sqrt n},
\end{equation}
with $c_0=1/\sqrt3$ by \eqref{eq:phasekhintchine}. H\"older's inequality in \eqref{eq:headtail} yields
\[
 U\le \left(\sum_{i<m}x_i^2\right)^{1/2}\left(\sum_{i<m}\frac1i\right)^{\beta/2}\le\norm{x}_2(1+\log m)^{\beta/2}.
\]
Take the probability measure $\eta=(1-\delta)\delta_y+\delta\delta_{y^{(0)}}$, where $\delta_z$ denotes the Dirac measure at $z$. If $U>\delta$, then
\begin{equation}\label{eq:headcontrol}
 \sqrt{\beta nH}\int A_n(x,z)\,d\eta(z)
 \ge\frac{c_0\sqrt\beta H^{1/2}\delta^2}{(1+\log m)^{\beta/2}}\longrightarrow\infty.
\end{equation}
Indeed, $\delta^2=H^{-1/4}$ and $Q^2\le\log m\le Q^2+\log 2$, so the right-hand side of \eqref{eq:headcontrol} is bounded below, for large $H$, by a positive constant depending only on $s$ times $H^{1/4}/(\log H)^{2\beta}$. If $U\le\delta$, \eqref{eq:lowerdense} bounds the same integral from below by $(1-\delta)A_n(x,y)$. Taking the infimum over the unit sphere, rearrangement invariance of the averages, \eqref{eq:lowerdense}, and \eqref{eq:headcontrol} give
\begin{equation}\label{eq:lowerorbitlimit}
 \liminf_{n\to\infty}\sqrt{\beta nH}\inf_{\norm{x}_{2,s}=1}\int A_n(x,z)\,d\eta(z)\ge\gamma_{\K},
\end{equation}
since the estimates hold for every $\varepsilon\in(0,\gamma_{\K})$. Equations \eqref{eq:common}, \eqref{eq:orbitexact}, \eqref{eq:upperorbitlimit} and \eqref{eq:lowerorbitlimit} now show that
\begin{equation}\label{eq:exactpi}
 \lim_{n\to\infty}\frac{\sqrt{\beta n\log\log n}}{\pi_1(\id_{\ell_{2,s}^n})}=\gamma_{\K}.
\end{equation}
By \eqref{eq:identity} and \eqref{eq:gammas}, this proves Theorem~\ref{thm:exact}, since $K_s^{\K}=\gamma_{\K}/\sqrt\beta$.

\section{Uniform estimates}
Let $1<s<2$, $t=s/(s-1)>2$ and $\beta=2/s-1=1-2/t\in(0,1)$. The estimates in this section have numerical constants. The endpoint cases $s=1,2$ follow by continuity for each fixed dimension. For $L\in\N$ and $n=2^L-1$, let
\begin{equation}\label{eq:blocks}
 I_k=\{2^k,\dots,2^{k+1}-1\}\quad (0\le k<L),
 \qquad H_L=\sum_{j=1}^L\frac1j.
\end{equation}
For $u_1\ge\cdots\ge u_n\ge0$, define
\[
 A_k(u)=2^{k/2}u_{2^{k+1}-1},\qquad B_k(u)=2^{k/2}u_{2^k}\quad(0\le k<L).
\]
\begin{lemma}\label{lem:dyad}
For every nonnegative non-increasing $u\in\R^n$,
\begin{align}
 \norm{u}_{2,s}&\asymp\left(\sum_{k=0}^{L-1}A_k(u)^s\right)^{1/s}
 \asymp\left(\sum_{k=0}^{L-1}B_k(u)^s\right)^{1/s},\label{eq:dyadprimal}\\
 \norm{u}_*&\asymp\left(\sum_{k=0}^{L-1}B_k(u)^t\right)^{1/t}.\label{eq:dyaddual}
\end{align}
\end{lemma}
\begin{proof}
Since $s/2-1<0$, for $i\in I_k$ one has
\[
 2^{(k+1)(s/2-1)}\le i^{s/2-1}\le 2^{k(s/2-1)}.
\]
Using $|I_k|=2^k$ and the monotonicity of $u$, we obtain
\[
 2^{s/2-1}A_k(u)^s\le\sum_{i\in I_k}i^{s/2-1}u_i^s\le B_k(u)^s.
\]
Also $B_0(u)=A_0(u)$ and $B_k(u)\le\sqrt2 A_{k-1}(u)$ for $1\le k<L$. In particular,
\begin{align*}
 2^{s/2-1}\sum_{k=0}^{L-1}A_k(u)^s
 &\le \norm{u}_{2,s}^s\le\sum_{k=0}^{L-1} B_k(u)^s\\
 &\le (1+2^{s/2})\sum_{k=0}^{L-1}A_k(u)^s.
\end{align*}
Since $A_k(u)\le B_k(u)$, we obtain \eqref{eq:dyadprimal}.

By the rearrangement inequality and H\"older's inequality, for any $x\in\K^n$,
\begin{align*}
 |\langle x,u\rangle|&\le\sum_{i=1}^n x_i^*u_i
 \le\sum_{k=0}^{L-1} B_k(x^*)B_k(u)\\
 &\le\left(\sum_{k=0}^{L-1}B_k(x^*)^s\right)^{1/s}
       \left(\sum_{k=0}^{L-1}B_k(u)^t\right)^{1/t}.
\end{align*}
Equation \eqref{eq:dyadprimal} gives the upper bound in \eqref{eq:dyaddual}.
For the lower bound write $d_k=A_k(u)$ and choose $c_k\ge0$ such that $\sum_k c_k^s=1$ and $\sum_k c_kd_k=(\sum_k d_k^t)^{1/t}$. Put
\[
 v_k=\max_{k\le j<L}2^{-(j-k)/2}c_j,
 \qquad x_i=2^{-k/2}v_k\quad(i\in I_k).
\]
The sequence $x$ is non-increasing. Extending $c_j$ by zero for $j\ge L$, Minkowski's inequality gives
\[
 \norm{v}_s\le\sum_{m=0}^{\infty}2^{-m/2}\left(\sum_{k=0}^{L-1}c_{k+m}^s\right)^{1/s}
 \le (1-2^{-1/2})^{-1}.
\]
Since $B_k(x)=v_k\ge c_k$, the block estimates above give $\norm{x}_{2,s}\le\norm{v}_s$ and $\langle x,u\rangle\ge\sum_k c_kd_k=\norm{d}_t$. Hence $\norm{u}_*\ge(1-2^{-1/2})\norm{d}_t$. The inequalities $B_0(u)=d_0$ and $B_k(u)\le\sqrt2 d_{k-1}$ yield
\begin{equation}\label{eq:dualnumeric}
 \left(\sum_k B_k(u)^t\right)^{1/t}\le\frac{(1+2^{t/2})^{1/t}}{1-2^{-1/2}}\norm{u}_*\le\frac{\sqrt3}{1-2^{-1/2}}\norm{u}_*<6\norm{u}_*.
\end{equation}
All constants are numerical, since $2^{s/2-1}\ge2^{-1/2}$, $1+2^{s/2}\le3$ and $t>2$.
\end{proof}

Let $n=2^L-1$, $L\ge2$, and $X=\ell_{2,s}^n$.
\subsection*{The upper bound for the absolutely 1-summing norm}
The choice of the dyadic weights is governed by the weighted $\ell_2$ norm with weights $(k+1)^\beta$. We take the squares of the weights to be successive differences of the power function:
\[
 d_l=\big((L-l)^\beta-(L-l-1)^\beta\big)^{1/2},\qquad
 y_i=2^{-l/2}d_l\quad(i\in I_l),\quad 0\le l<L.
\]
With $0^\beta=0$, these weights satisfy
\begin{equation}\label{eq:tailuniform}
 \sum_{l=L-1-k}^{L-1}d_l^2=(k+1)^\beta\qquad(0\le k<L).
\end{equation} For $j\ge2$, since $\beta-1=-2/t$,
\[
 j^\beta-(j-1)^\beta=\int_{j-1}^{j}\beta u^{\beta-1}\,du
 \le\beta(j-1)^{-2/t}.
\]
The term corresponding to $j=1$ is $1$, and therefore
\begin{align}
 \sum_{l=0}^{L-1}d_l^t
 &=\sum_{j=1}^{L}\big(j^\beta-(j-1)^\beta\big)^{t/2}\notag\\
 &\le 1+\beta^{t/2}\sum_{j=2}^L(j-1)^{-1}\le1+\beta^{t/2}H_L.\label{eq:dnorm}
\end{align}
For any $u\in\K^n$, the rearrangement inequality followed by H\"older gives
\begin{equation}\label{eq:directdual}
 \norm{u}_*\le\left(\sum_{i=1}^n i^{t/2-1}(u_i^*)^t\right)^{1/t}.
\end{equation}
Indeed, the factors in the sum $\sum_i x_i^*u_i^*$ are $i^{1/2-1/s}x_i^*$ and $i^{1/s-1/2}u_i^*$, and $t(1/s-1/2)=t/2-1$. Since $t>2$, the weights $i^{t/2-1}$ are increasing; pairing them with the decreasing rearrangement minimizes the weighted sum. Consequently
\begin{align}\label{eq:ynormuniform}
 \norm{y}_*^t&\le\sum_{i=1}^n i^{t/2-1}|y_i|^t\notag\\
 &\le\sum_{l=0}^{L-1}2^{(l+1)(t/2-1)}2^l 2^{-lt/2}d_l^t\notag\\
 &=2^{t/2-1}\sum_{l=0}^{L-1}d_l^t\le 2^{t/2-1}(1+\beta^{t/2}H_L).
\end{align}

Let $x\in X$. By invariance we may suppose $x$ is nonnegative and non-increasing. Write $a_k=A_k(x)$. For $i\in I_k$ and $j\in I_l$ define
\[
 w_{ij}=\begin{cases}
 2^{-k-l}a_k^2d_l^2,& k+l\ge L-1,\\
 0,& k+l<L-1.
 \end{cases}
\]
Then $w_{ij}\le|x_i y_j|^2$. By \eqref{eq:momentmean} and $|I_k|=2^k$, the mean $M$ of $Z(\sigma)=\sum_{i=1}^n w_{i,\sigma(i)}$ is
\begin{equation}\label{eq:mean}
 M=\frac1n\sum_{k+l\ge L-1}a_k^2d_l^2
 =\frac1n\sum_{k=0}^{L-1}a_k^2\sum_{l=L-1-k}^{L-1}d_l^2.
\end{equation}
For a nonzero entry in a retained block, $2^{k+l}\ge2^{L-1}>n/2$. The contribution of that entire block to $M$ is $2^{k+l}w_{ij}/n>w_{ij}/2$. It follows that every entry satisfies $w_{ij}\le2M$. Lemma~\ref{lem:moment} and \eqref{eq:F} yield
\begin{equation}\label{eq:FLower}
 F_n(x,y)\ge\frac{1}{2\sqrt n}\left(\sum_{k=0}^{L-1}a_k^2\sum_{l=L-1-k}^{L-1}d_l^2\right)^{1/2}.
\end{equation}
By \eqref{eq:tailuniform}, the right-hand side of \eqref{eq:FLower} is
\begin{equation}\label{eq:FLower2}
 F_n(x,y)\ge\frac1{2\sqrt n}\left(\sum_{k=0}^{L-1}a_k^2(k+1)^\beta\right)^{1/2}.
\end{equation}
H\"older's inequality with exponents $2/s$ and $2/(2-s)$ gives
\[
 \left(\sum_{k=0}^{L-1}a_k^s\right)^{1/s}
 \le H_L^{\beta/2}\left(\sum_{k=0}^{L-1}a_k^2(k+1)^\beta\right)^{1/2}.
\]
By \eqref{eq:FLower2} and the bounds in the proof of Lemma~\ref{lem:dyad},
\[
 F_n(x,y)\ge\frac{\norm{x}_{2,s}}{2(1+\sqrt2)\sqrt n\,H_L^{\beta/2}}.
\]
Here $(1+2^{s/2})^{1/s}\le1+\sqrt2$ by the monotonicity of finite-dimensional $\ell_p$ norms. On the other hand, \eqref{eq:ynormuniform}, the inequality $(a+b)^{2/t}\le a^{2/t}+b^{2/t}$ for $a,b\ge0$, and $H_L^\beta\le1-\beta+\beta H_L$ give
\begin{align*}
 H_L^\beta\norm{y}_*^2
 &\le 2H_L^\beta(1+\beta^{t/2}H_L)^{2/t}\\
 &\le 2(H_L^\beta+\beta H_L^{\beta+2/t})\\
 &\le 2(1+2\beta H_L)\le4(1+\beta H_L).
\end{align*}
Applying \eqref{eq:upperpi} to $y/\norm{y}_*$ with $c=1/\sqrt3$ yields
\begin{equation}\label{eq:PIupper}
 \pi_1(\id_X)\le4\sqrt3(1+\sqrt2)\sqrt{n(1+\beta H_L)}.
\end{equation}

\subsection*{The lower bound for the absolutely 1-summing norm}
Define a nonnegative non-increasing vector $x$ by
\begin{equation}\label{eq:testx}
 x_i=2^{-k/2}a_k\quad (i\in I_k),\qquad a_k=(k+1)^{-1/s}.
\end{equation}
Since $A_k(x)=B_k(x)=a_k$, \eqref{eq:dyadprimal} implies
\begin{equation}\label{eq:testnorm}
 2^{-1/2}H_L^{1/s}\le\norm{x}_{2,s}\le H_L^{1/s}.
\end{equation}
Fix $y\in B_{X^*}$. We may suppose $y$ nonnegative and non-increasing, because both the dual norm and $F_n$ are invariant under rearrangements of the absolute values. Set $b_l=B_l(y)$. Then $y_j\le2^{-l/2}b_l$ for $j\in I_l$, and \eqref{eq:dualnumeric} gives
\begin{equation}\label{eq:bdual}
 \left(\sum_{l=0}^{L-1}b_l^t\right)^{1/t}\le6.
\end{equation}
For $\sigma\in S_n$ and $0\le k,l<L$, write
\[
 N_{kl}(\sigma)=|\{i\in I_k:\sigma(i)\in I_l\}|.
\]
Since $\mu_n\{\sigma:\sigma(i)=j\}=1/n$,
\begin{equation}\label{eq:Nmean}
 \int_{S_n}N_{kl}\,d\mu_n=\frac{2^{k+l}}n.
\end{equation}
By \eqref{eq:testx}, for each permutation,
\begin{equation}\label{eq:XYupper}
 \sum_{i=1}^n|x_i y_{\sigma(i)}|^2\le\sum_{k,l=0}^{L-1}2^{-k-l}a_k^2b_l^2 N_{kl}(\sigma).
\end{equation}
Put $r_l=L-1-l$. Split the indices on the right-hand side of \eqref{eq:XYupper} into
\[
 \mathcal D=\{(k,l):0\le k,l<L,\ 2k\ge r_l\},\quad
 \mathcal E=\{(k,l):0\le k,l<L,\ 2k<r_l\}.
\]
Put
\begin{align*}
 U&=\int_{S_n}\left(\sum_{(k,l)\in\mathcal D}2^{-k-l}a_k^2 b_l^2 N_{kl}(\sigma)\right)^{1/2}d\mu_n(\sigma),\\
 V&=\int_{S_n}\left(\sum_{(k,l)\in\mathcal E}2^{-k-l}a_k^2 b_l^2 N_{kl}(\sigma)\right)^{1/2}d\mu_n(\sigma).
\end{align*} The inequality $\sqrt{a+b}\le\sqrt a+\sqrt b$ for $a,b\ge0$ gives
\begin{equation}\label{eq:UV}
 F_n(x,y)\le U+V.
\end{equation}
Jensen's inequality and \eqref{eq:Nmean} give
\[
 U\le\frac1{\sqrt n}\left(\sum_{l=0}^{L-1}b_l^2\sum_{k\ge r_l/2}^{L-1}(k+1)^{-2/s}\right)^{1/2},
\]
where the inner sum is over integers $k$. For every integer $m\ge0$,
\[
 \sum_{k=m}^{\infty}(k+1)^{-2/s}
 \le(m+1)^{-2/s}+\int_{m+1}^{\infty}u^{-2/s}\,du
 \le(1+\beta^{-1})(m+1)^{-\beta}\le\frac{2}{\beta}(m+1)^{-\beta}.
\]
For $m=\lceil r_l/2\rceil$ we have $m+1\ge(r_l+1)/2$. As $1-2/s<0$, the last estimate gives
\begin{equation}\label{eq:Ubound}
 U\le\frac2{\sqrt{n\beta}}\left(\sum_{l=0}^{L-1}b_l^2(r_l+1)^{-\beta}\right)^{1/2}.
\end{equation}
By H\"older's inequality and \eqref{eq:bdual}, since $\beta=1-2/t$,
\[
 \sum_{l=0}^{L-1}b_l^2(r_l+1)^{-\beta}
 \le\norm{b}_t^2 H_L^\beta\le36H_L^\beta.
\]
Thus \eqref{eq:Ubound} yields
\begin{equation}\label{eq:Ufinal}
 U\le\frac{12H_L^{\beta/2}}{\sqrt{n\beta}}.
\end{equation}

For the sum over $\mathcal E$, use the subadditivity of the square root and the fact that $\sqrt N\le N$ for every nonnegative integer $N$. By \eqref{eq:Nmean},
\begin{align*}
 V&\le\sum_{(k,l)\in\mathcal E}2^{-(k+l)/2}a_k b_l\int_{S_n}\sqrt{N_{kl}}\,d\mu_n\\
 &\le\frac1n\sum_{(k,l)\in\mathcal E}2^{(k+l)/2}a_k b_l.
\end{align*}
Since $2^{L-1}\le n$, we have $2^{(k+l)/2}/n\le n^{-1/2}2^{(k-r_l)/2}$. By \eqref{eq:ycoord}, $b_l\le2$. Using $a_k\le1$ and extending the sums to all nonnegative integers, we obtain
\begin{align*}
 V&\le\frac{2}{\sqrt n}\sum_{k=0}^{\infty}\sum_{r=2k+1}^{\infty}2^{(k-r)/2}\\
 &=\frac{2\cdot2^{-1/2}}{(1-2^{-1/2})^2\sqrt n}<\frac{17}{\sqrt n}.
\end{align*}
Since $H_L\ge1$ and $0<\beta<1$, \eqref{eq:UV} and \eqref{eq:Ufinal} give
\begin{equation}\label{eq:Fupper}
 \sup_{y\in B_{X^*}}F_n(x,y)\le\frac{29H_L^{\beta/2}}{\sqrt{n\beta}}.
\end{equation}
By \eqref{eq:lowerpi}, \eqref{eq:testnorm} and \eqref{eq:Fupper},
\[
 \pi_1(\id_X)\ge\frac{\sqrt{n\beta H_L}}{29\sqrt2}.
\]
Lemma~\ref{lem:identity} and the Kadets--Snobar theorem also give $\pi_1(\id_X)\ge\sqrt n$. Adding the two squared bounds in the form $n\le\pi_1(\id_X)^2$ and $n\beta H_L\le 2\cdot29^2\pi_1(\id_X)^2$, we obtain
\begin{equation}\label{eq:PIlower}
 \pi_1(\id_X)\ge\frac{\sqrt{n(1+\beta H_L)}}{\sqrt{1+2\cdot29^2}}\ge\frac{\sqrt{n(1+\beta H_L)}}{42}.
\end{equation}
By \eqref{eq:PIupper} and \eqref{eq:PIlower},
\begin{equation}\label{eq:piorder}
 \frac1{42}\sqrt{n(1+\beta H_L)}\le\pi_1(\id_X)\le 4\sqrt3(1+\sqrt2)\sqrt{n(1+\beta H_L)}.
\end{equation}

\subsection*{Passage to arbitrary dimensions}
For $1\le m\le n$, zero extension $J\colon\ell_{2,s}^m\to\ell_{2,s}^n$ is isometric. Applying \eqref{eq:pi} to the images of a finite family under $J$ and restricting the dual functionals to $J(\ell_{2,s}^m)$ gives
\begin{equation}\label{eq:monotonepi}
 \pi_1(\id_{\ell_{2,s}^m})\le\pi_1(\id_{\ell_{2,s}^n}).
\end{equation}
Given $n\ge3$, choose $L\ge2$ so that $m=2^L-1\le n<M=2^{L+1}-1$. Then $n\le2m$, $M=2m+1\le2n+1\le7n/3$, and \eqref{eq:monotonepi} gives
\[
 \pi_1(\id_{\ell_{2,s}^{m}})\le\pi_1(\id_{\ell_{2,s}^{n}})\le\pi_1(\id_{\ell_{2,s}^{M}}).
\]
The integral comparison $\log(L+1)\le H_L\le1+\log L$ and the inequalities
\[
 e+\log n\le e+(L+1)\log2\le2(L+1),\qquad
 e+\log n\ge e+(L-1)\log2\ge\sqrt{e(L+1)}
\]
for $L\ge2$ give $h(n)\le2H_L$ and $H_{L+1}\le2h(n)$. For the last elementary inequality, the difference $(e+(L-1)\log2)^2-e(L+1)$ is positive and increasing for $L\ge2$.
Applying \eqref{eq:piorder} at $m$ and $M$, with $m\ge n/2$ and $M\le7n/3$, gives
\begin{align*}
 \pi_1(\id_{\ell_{2,s}^n})&\ge \frac1{84}\sqrt{n(1+\beta h(n))},\\
 \pi_1(\id_{\ell_{2,s}^n})&\le4\sqrt3(1+\sqrt2)\sqrt{14/3}\sqrt{n(1+\beta h(n))}\\
 &\le40\sqrt{n(1+\beta h(n))}.
\end{align*}
By \eqref{eq:identity}, these are \eqref{eq:uniformbounds} with $c=1/40$ and $C=84$ for $n\ge3$. For $n=1,2$, $h(n)\le2$ and $1\le\lambda(\ell_{2,s}^n)\le\sqrt n$, which imply the same bounds.

For fixed $n$, compactness of the Euclidean unit sphere and continuity of the norm in \eqref{eq:lorentz} imply that $\norm{\cdot}_{2,s}/\norm{\cdot}_{2,s_0}\to1$ uniformly on nonzero vectors as $s\to s_0\in\{1,2\}$. The identity $T_s\colon\ell_{2,s}^n\to\ell_{2,s_0}^n$ therefore satisfies $\norm{T_s}\norm{T_s^{-1}}\to1$. The ideal property gives
\[
 \frac{\pi_1(\id_{\ell_{2,s_0}^n})}{\norm{T_s}\norm{T_s^{-1}}}
 \le\pi_1(\id_{\ell_{2,s}^n})
 \le\norm{T_s}\norm{T_s^{-1}}\pi_1(\id_{\ell_{2,s_0}^n}).
\]
By Lemma~\ref{lem:identity}, $\lambda(\ell_{2,s}^n)\to\lambda(\ell_{2,s_0}^n)$, and the uniform bounds extend to both endpoints. This completes the proof of Theorem~\ref{thm:main}.

\section*{Declarations}
\noindent\textbf{Funding.} D. N\'u\~nez-Alarc\'on, D. M. Pellegrino and D. M. Serrano-Rodr\'iguez are supported by the Conselho Nacional de Desenvolvimento Cient\'ifico e Tecnol\'ogico (CNPq, Brazil) through Grants No.~406457/2023-9 (CNPq/MCTI Call No.~10/2023) and No.~403964/2024-5 (MCTI/CNPq Call No.~16/2024). D. M. Pellegrino is also supported by CNPq Grant No.~305807/2025-0.

\medskip
\noindent\textbf{Competing interests.} The authors declare that they have no competing interests.

\medskip
\noindent\textbf{Data availability.} No datasets were generated or analysed in this study.

\medskip
\noindent\textbf{Declaration of generative AI and AI-assisted technologies.} OpenAI's ChatGPT was used in the exploration of the problem, the development of mathematical arguments, bibliographic searches and the preparation and revision of the manuscript. The authors assume responsibility for the mathematical content and the final version of the manuscript.


\end{document}